                     \def\version{May 21st, 2012}                   %
\documentclass[reqno,11pt]{amsart} 
\usepackage[notref,notcite]{showkeys}
\usepackage{srcltx}

\usepackage{amsmath, amsthm, a4, latexsym, amssymb}
\usepackage{graphicx}
\makeatletter\@addtoreset{equation}{section}\makeatother

 \newfont{\bfit}{cmbxti10 scaled 1200}

 \newcommand{\e}{{\rm e} }

 \newcommand{\N}{\mathbb{N}}
 \newcommand{\Z}{\mathbb{Z}}

 \newcommand{\E}{\mathbb{E}}
 
 \newcommand{\tL}{\tilde L}
 \newcommand{\tU}{\tilde U}
 \newcommand{\tlambda}{\tilde\lambda}
 \renewcommand{\P}{\mathbb{P}}
 \def\1{{\mathchoice {1\mskip-4mu\mathrm l} 
{1\mskip-4mu\mathrm l}
{1\mskip-4.5mu\mathrm l} {1\mskip-5mu\mathrm l}}}

 \newcommand{\BbbP}{\mathbb P}
 \newcommand{\BbbE}{\mathbb E}
 \newcommand{\quenchedP}{P}
 \newcommand{\annealedP}{\BbbP}
 \newcommand{\quenchedE}{E}
 \newcommand{\annealedE}{\BbbE}
 \newcommand{\ignore}[1]{}

\renewcommand{\subsection}{\secdef \subsct\sbsect}
\newcommand{\subsct}[2][default]{\refstepcounter{subsection}
\vspace{0.15cm}
{\flushleft\bf \arabic{section}.\arabic{subsection}~\bf #1  }
\nopagebreak\nopagebreak}
\newcommand{\sbsect}[1]{\vspace{0.1cm}\noindent
{\bf #1}\vspace{0.1cm}}

{\nopagebreak {\hfill\rule{2mm}{2mm}}\\ }

\newtheorem{theorem}{Theorem}[section]
\newtheorem{lemma}[theorem]{Lemma}

\newtheorem{prop}[theorem]{Proposition}

\newtheorem{remark}[theorem]{Remark}

\newtheoremstyle{thm}{1.5ex}{1.5ex}{\itshape\rmfamily}{}
{\bfseries\rmfamily}{}{2ex}{}

\newtheoremstyle{rem}{1.3ex}{1.3ex}{\rmfamily}{}
{\itshape\rmfamily}{}{1.5ex}{}
\theoremstyle{rem}
\refstepcounter{subsubsection}

\def\thebibliography#1{\section*{References}
  \list%
  {\arabic{enumi}.}%                          {\star}{\star}{\star} style of reference number {\star}{\star}{\star}
    {\settowidth\labelwidth{[#1]}\leftmargin\labelwidth
    \advance\leftmargin\labelsep
    \parsep0pt\itemsep0pt
    \usecounter{enumi}}
    \def\newblock{\hskip .11em plus .33em minus .07em}
    \sloppy                   % \clubpenalty4000\widowpenalty4000
    \sfcode`\.=1000\relax}

\begin{document}
%%%%%%%%%%%%%%%%%%%%%%%%%%%%%%%%%%%%%%%%%%%%%%%
\title[RWRC with no speed]
{\large On the speed of Random Walks among Random Conductances}
\author[Noam berger and Michele Salvi]{}
\maketitle
\thispagestyle{empty}
\vspace{-0.5cm}

\centerline{\sc By Noam Berger$ $\footnote{Einstein Institute of Mathematics, The Hebrew University, Givat Ram, Jerusalem 91904, Israel. \tt{berger@math.huji.ac.il}} and Michele Salvi\footnote{Institute for Mathematics, TU Berlin, Str.~des 17.~Juni 136, 10623 Berlin, Germany. {\tt salvi@math.tu-berlin.de}}}
\renewcommand{\thefootnote}{}


\vspace{-0.5cm}
\centerline{\textit{Hebrew University of Jerusalem and TU Berlin}}
\vspace{0.2cm}

\begin{center}
\version
\end{center}

\begin{abstract}
We consider random walk among random conductances where the conductance environment is shift invariant and ergodic. 
We study which moment conditions of the conductances guarantee speed zero of the random walk. 
We show that if there exists $\alpha>1$ such that $E[\log^\alpha \omega_e]<\infty$, then the random walk has speed zero. 
On the other hand, for each $\alpha<1$ we provide examples of random walks with non-zero speed and random walks for which the limiting speed does not exist that have $E[\log^\alpha \omega_e]<\infty$.
%We also provide examples of a random walk with non-zero speed and one for which the limiting speed does not exist
%that have $E[\log^\alpha \omega_e]<\infty$ for each $\alpha<1$.
\end{abstract}

\section{Introduction}
\subsection{Definitions and notations}
We discuss the following two-dimensional model of motion in random medium. 
Let $\E^2$ be the set of nearest neighbor edges in the graph $\Z^2$. 
We may also write an edge as an unordered pair $\{x,y\}$ of neighbors in $\Z^2$. 
We consider the measure space $\Omega=(0,\infty)^{\E^2}$. 
For $\omega\in\Omega$ and an edge $e$, we call $\omega_e$ the {\em weight} or {\em conductance} of the edge $e$.
We let $P$ be a measure on $\Omega$ which satisfies the following two conditions:

\begin{itemize}
\item $P$ is invariant and ergodic w.r.t. the group of spatial moves in $\Z^2$.
\item The marginal distribution of $\omega_e$ is the same for all choices of the edge $e$, i.e.~vertical and horizontal edges have the same distribution.
\end{itemize}
Note that this is weaker than invariance w.r.t. rotations. The second condition can be weakened significantly, but for simplicity we keep it as is. 

For a given $\omega\in\Omega$ and $z\in\Z^2$ we define the law $\quenchedP^\omega_z$ as follows. 
$\quenchedP^\omega_z$ is the law of a Markov chain on $\big(\Z^2\big)^\N$ with 
$\quenchedP^\omega_z(X_0=z)=1$ and 
$$
\quenchedP^\omega_z(X_{n+1}=y|X_n=x)
=\frac{\omega(x,y)}{\sum_{w\sim x}\omega(x,w)}
$$ 
for any neighbouring points $x,y\in\Z^2$, and where the sum in the denominator is carried over every neighbour $w$ of $x$.

For any event $A\subseteq(\Z^2)^\N$, we define the annealed law to be 
\[
\annealedP_z(A) = \int_\Omega \quenchedP^\omega_z(A) P(d\omega).
\]

If $z$ is the origin, we may omit the subscript. We use $\annealedE_z$ and $\quenchedE^\omega_z$ as the 
expectations w.r.t. the distributions $\annealedP_z$ and $\quenchedP^\omega_z$.

This is a well-known process, called the random walk among random conductances (RWRC).
The RWRC has been studied extensively in the past decades, see e.g.~\cite{Biskup} and the references therein. It is related to many other important models in Statistical physics, for example the gradient fields (see e.g.~\cite{Biskup-Spohn}) and reinforced random walk (see e.g.~\cite{Merkl-Rolles, Angel-Crawford-Kozma, Sabot-Tarres}).

Traditionally, RWRC is studied in connection to problems such as recurrence and transience \cite{BergerLRP, Angel-Crawford-Kozma, Sabot-Tarres}, central limit theorems (see \cite{Biskup-Prescott, Mathieu, Andres-Barlow-Deuschel-Habmbley} and many more) and heat kernel estimates (see e.g.~\cite{Delmotte, Barlow-Bass, BBHK08}).
However, in the present paper we ask a slightly different question, that of the limiting velocity of the random walk.

\subsection{Main question}
It is easy to prove and well known that if the conductances are bounded then the speed is zero, i.e.~$\annealedP(\lim(X_n/n)=0)=1$. It is also well known that if the conductances are i.i.d.~the same is true, regardless of the distribution of the single conductances. We thus wish to understand which conditions force the speed of the RWRC to be zero. Based on the examples above, it seems that two types of criteria are involved. The first is moment conditions that control the size of the conductances, and the second is mixing conditions saying that if the environment mixes fast enough then the speed is zero.

In this paper we only consider the first type, and show that the sharp condition is that the logarithm of the conductances has high enough moments.

Our main result is as follows.
\begin{theorem}\label{thm:main} Let $e$ be an edge in $\E^2$.
\begin{itemize}
\item[(i)]\label{item:carne} If there exists $\alpha>1$ such that 
\begin{equation}\label{eq:logmomalpha}
E[\log^\alpha\omega_e]<\infty,
\end{equation}
 then
\[
\annealedP\left(\lim_{n\to\infty}\frac{X_n}{n}=0\right)=1.
\]
\item[(ii)]\label{item:trees} For every $\alpha<1$ there exists a distribution $P$ on environments such that 
$E[\log^\alpha\omega_e]<\infty$, but
\[
\annealedP\left(\lim_{n\to\infty}\frac{X_n}{n}=0\right)=0.
\]
Furthermore, in this case it is possible to choose $P$ so that either
$$
\annealedP\left(\Big\|\lim_{n\to\infty}\frac{X_n}{n}\Big\|_\infty>0\right)=1
$$
or
$$
\annealedP\left(\lim_{n\to\infty}\frac{X_n}{n} \mbox{ does not exist}\right)=1.
$$
\end{itemize}
\end{theorem}

\begin{remark}
Our proofs deal with conductances bounded away from zero, but would work in the same way including the possibility of zero conductances. Note also that the choice of dimension $2$ has been made in order to have easier and more intuitive proofs. We are confident that the same results can be proven with the very same techniques in higher dimensions, with critical $\alpha$ equal to $d-1$.
\end{remark}
\begin{remark}
Our counter examples are not uniformly elliptic (i.e.~the transition probabilities are not bounded away from zero). A natural question is whether it is possible to construct similar examples such that the transition probabilities are uniformly elliptic (cf. \cite{heil}).
\end{remark}

In Section \ref{sec:carne} we show Part (i) of Theorem \ref{thm:main}, which ends up being a simple application of the Varopoulos-Carne Theorem. In Sections \ref{sec:trees} and \ref{sec:env} we show Part (ii) of Theorem \ref{thm:main}. The construction builds upon the example constructed by Bramson, Zeitouni and Zerner in \cite{BZZ}.

\section{Moment conditions for speed zero}\label{sec:carne}
In this section we prove Part (i) of Theorem \ref{thm:main}.
%\begin{prop}\label{prop:carne}
%Let $\alpha>1$ and assume that $P$ is such that 
%\begin{equation}\label{eq:logmomalpha}
%E\big[\log^\alpha\omega_e\big]<\infty.
%\end{equation}
%Then
%\[
%\annealedP\Big(\lim_{n\to\infty}\frac{X_n}{n}=0\Big)=1.
%\]
%\end{prop}

In order to prove it, we will use the well known Varopoulos-Carne bound. For proof see, e.g., \cite{carne}.

\begin{lemma}[Varopoulos-Carne]\label{lem:carne}
Let $L$ be an irreducible Markov transition kernel with reversible measure $\pi$. For states $x$ and $y$, denote
$d(x,y)=\min\{n:L^n(x,y)>0\}$.
Then for every $x$, $y$ and $n$,
\begin{equation}\label{eq:carne}
L^n(x,y)\leq
2\sqrt{\tfrac{\pi(y)}{\pi(x)}}\cdot e^{-\frac{d(x,y)^2}{2n}}.
\end{equation}
\end{lemma}

\begin{proof}[Proof of Part (i) of Theorem \ref{thm:main}]
The measure $\pi$ on $\Z^2$, defined by $\pi(x)=\sum_{y\sim x} \omega_{\{x,y\}}$ is a reversible measure for our random walk. 
As in \eqref{eq:logmomalpha},
let
\[
D = E\left[\log^\alpha\omega_e\right]<\infty.
\]

For $n\in \N$, consider the points $x\in\Z^2$ such that $||x||_\infty=n$, and call $E_n$ the set of edges having at least one end in these points. 
Note that $|E_n|=24n$.

Then by Markov's inequality, for every $n\in\N$ and $K>0$ we get 
\begin{eqnarray*}
P\big(
\exists e\in E_n \mbox{ s.t. } \omega_e>\tfrac K 4
\big)
\leq 24n \frac{ D}{\log^\alpha (\tfrac K 4)}.
\end{eqnarray*}

In particular, if $K=e^{n^\beta}$ with $1/\alpha<\beta<1$, then 
\begin{eqnarray*}
P\left(
\exists e\in E_n \mbox{ s.t. } \omega_e>\tfrac K 4
\right)
\leq C n^{1-\alpha\beta},
\end{eqnarray*}
for some constant $C>0$.

Observe that $1-\alpha\beta<0$. Therefore, by the Borel-Cantelli lemma, for an integer $\kappa>(\alpha\beta -1)^{-1}$, 
a.s. for all $n$ large enough and every edge $e\in E_{n^\kappa}$, we have
\[
\omega_e\leq\tfrac 1 4 e^{n^{\kappa\beta}}.
\]

Therefore, for every $x$ s.t. $\|x\|_\infty=n^\kappa$, we have that $\pi(x)\leq e^{n^{\kappa\beta}}.$

Now fix $M\in\N$ and assume that $M$ is large. For every $n$ large enough,

\begin{align*}
\quenchedP^\omega\big(\|X_{Mn^\kappa}\|_\infty>n^\kappa\big)
&\leq \quenchedP^\omega\big(\exists{k\leq Mn^\kappa}:\,\|X_{Mn^\kappa}\|_\infty=n^\kappa\big)\\
&\leq \sum_{k=1}^{Mn^\kappa}\sum_{x:\,\|x\|_\infty=n^\kappa}\quenchedP^\omega(X_k=x)\\
&\leq \sum_{k=1}^{Mn^\kappa}\sum_{x:\,\|x\|_\infty=n^\kappa} 2\sqrt{\tfrac{\pi(x)}{\pi(0)}}
e^{-\tfrac{n^{2\kappa}}{2k}}\\
&\leq C'\pi(0)^{-1/2} \exp\big\{{\tfrac {n^{\kappa\beta}}{2}-\tfrac{n^{\kappa}}{2M}}\big\},
\end{align*}
for some constant $C'>0$.

Therefore, again by Borel-Cantelli, almost surely for all $n$ large enough,
\[
\|X_{Mn^\kappa}\|_\infty \leq n^\kappa.
\]

From here we immediately get that almost surely
\[
\limsup_{n\to\infty}\frac{\|X_n\|_\infty}{n} \leq \frac{2}{M}
\]

and in fact, since $M$ is arbitrary, 

\[
\annealedP\Big(
\lim_{n\to\infty}\frac{X_n}{n} = 0
\Big)=1.
\]

\end{proof}

\section{Trees}\label{sec:trees}

In this section and in the next one we prove Part (ii) of Theorem \ref{thm:main}. The section is divided into two different subsections. In Subsection \ref{sec:nospeed} we create the structure for the random environment where, with probability one, the sequence $\big(\frac{X_n}{n}\big)$ does not converge, and in Subsection \ref{sec:speedpos} we create another example where with probability one the sequence $\big(\frac{X_n}{n}\big)$ converges to a speed which is not zero. In both cases $E[\log^{\alpha}\omega_e]<\infty$
for arbitrary $\alpha<1$. The example in Subsection \ref{sec:nospeed} is a direct application of the tree construction of Bramson, Zeitouni and Zerner \cite{BZZ}. For the construction in Subsection \ref{sec:speedpos}, we need to modify the tree of \cite{BZZ}. The construction is inspired by the construction in \cite{BZZ}, but we need to change quite a few details in order for the speed to converge.

In both cases, we adapt trees into environments for the random walk in the exact same fashion. This is done in Section \ref{sec:env}. Now, we give a short introduction with the necessary terms from \cite{BZZ}, and then, in Subsection \ref{sec:nospeed} and \ref{sec:speedpos}, we create the actual trees.

%We will restrict ourselves to the two dimensional case.
An \textit{ancestral function} is a (in our case random) function $a:x\in\Z^2\to a(x)\in \Z^2$ with the following properties:
\begin{itemize}
 \item $x$ and $a(x)$ are nearest neighbours;
 \item $a(a(x))\not =x$;
 \item the set of edges $F_a:=\{\{x,a(x)\}: x\in\Z^2\}$ is a forest (i.e.~the graph $(\Z^2,F_a)$ contains no cycles).
\end{itemize}
Every connected component of $F_a$ is an infinite tree. $a(x)$ can be seen as the parent of $x$ and we denote by $a^n(x)$ the $n$-th generation ancestor of $x$, for $n\geq 0$ (with the convention $a^0(x)=x$).

We also say that an ancestral function is directed if for some $i,j\in\{+1,-1\}$ and for every $x\in Z^2$, $a(x)-x\in\{(0,i),(j,0)\}$.

The length of the longest branch starting in $x$ (or the distance from $x$ of its farthest descendant, if one prefers the genealogical metaphore) is
\begin{equation}\label{eqn:bulk}
 h(x):=\sup\{n\geq 0: \exists y\in\Z^2 \mbox{ such that } a^n(y)=x\}.
\end{equation}

We are interested in the distribution of $h(0)$ in the case of a random translation invariant ancestral function.

Theorem 1 in \cite{BZZ} says that for any stationary ancestral function there exists a constant $c\geq 0$ such that
\begin{equation}\label{eqn:lower_bound_h}
 \liminf_{n\to\infty} n P(h(0)\geq n)\geq c.
\end{equation}

In the same article the authors show that this is in fact the best lower bound achievable. We give the 2-dimensional version of Theorem 2 in that paper:

\begin{theorem}[\cite{BZZ}, Theorem 2]\label{thm:BZZ}
 There exists a stationary directed ancestral function $(a(x))_{x\in\Z^2}$ that is polynomially mixing of order 1 and for which 
\begin{equation}\label{eqn:upper_bound_h}
 \limsup_{n\to\infty} n P(h(0)\geq n)<\infty.
\end{equation}
\end{theorem}

%As said, in the next section we will report the explicit construction of such an ancestral function. The tree obtained this way (hereafter BZZ-tree) will be the support for an opportunely designed random walk among random conductances that has almost surely no speed.
%In section ??? we will modify the example of Bramson, Zeitouni and Zermer in order to obtain a %new tree (Diagonal tree) with a tail behaviour of the function in (\ref{eqn:bulk}) similar to that of %the BZZ-tree. Nevertheless, the random walk among random conductances built on this new tree %will have an almost sure asymptotic speed which is different from zero.

We now describe the BZZ tree, as appearing in \cite{BZZ}.

\subsection{The BZZ tree}\label{sec:nospeed}

We provide now the construction of the ancestral function used in \cite{BZZ}, restricted to the 2 dimensional case. We will make use of the same notations as \cite{BZZ} with an additional tilde.

Let $\{e_1, e_2\}$ be the canonical basis of $\Z^2$, with $e_1$ parallel to the $x$-axis. Fix two constants $\tilde\theta$ and $\tilde n_0\in\N$ such that $2\sqrt 2\leq\tilde\theta\leq\tilde n_0^2$. For every $x\in\Z^2$ let $\tL(x)$ be i.i.d.~random variables with atomless distribution and satisfying
\begin{equation}
 \tilde  P(\tL(x)>t)= \frac {\tilde \theta} {t^2} \qquad \mbox{for } t\geq \tilde n_0.
\end{equation}

We define an \textit{umbrella} of intesity $t$ to be
\begin{equation}
 \tU_t=\bigcup_{i=1,2} \tU_{i,t}
\end{equation}
where
\begin{equation}
 \tU_{i,t}=\big\{y=(y_1,y_2)\in\Z^2:y_i=0,y_j\in(0,t], j\not = i\big\}
\end{equation}

are the sides of the umbrella. The strength of the umbrella is also defined to be equal to its intensity. 

For every $x\in\Z^2$ we will open the umbrella $x+\tU_{\tL(x)}$. Informally, one can think of the ancestral function as a drop of rain trying to fall towards the up-right direction of the plane and sliding on the sides of the umbrellas. Whenever two or more umbrellas overlap, the water will consider only the strongest of them and penetrate the perpendicular ones.

Formally, one defines for every $x\in\Z^2$ the strongest umbrella passing through that point perpendicular to direction $e_i$, for $i\in\{1,2\}$, as
\begin{equation}
 \tlambda_i(x)=\sup_{y\in\Z^2:\,x\in y+\tU_{i,\tL(y)}}\tL(y).
\end{equation}

Note that the $\sup$ is taken over a non-empty set and it is easy to show that $\tlambda_i(x)$ is also a.s. finite.

Since the distributions of the $\tL(x)$'s are atomless, the direction $I(x)\in\{1,2\}$ such that
$$
\tlambda_{I(x)}(x)=\min\{\tlambda_i(x), i=1,2\}
$$ 
is well defined.
The ancestral function we are looking for is
\begin{equation}
 \tilde a(x)=x+e_{I(x)}.
\end{equation}
The set of edges $\big\{\{x,\tilde a(x)\}, \,x\in\Z^2\big\}$ through which the drops of rain have flown forms a random forest (which can be shown to be in fact a random tree spanning the whole $\Z^2$). This is the ancestral function used to prove Theorem \ref{thm:BZZ}, and we will call the graph obtained with it the BZZ tree.

\subsection{The Diagonal tree}\label{sec:speedpos}

%\begin{figure}
%\begin{center}
%\includegraphics[scale=0.3]{Umbrellas}
% \caption{\textit{Both in the straight umbrellas case} $\rm{a)}$ \textit{and in the narrow umbrellas case} $\rm{b)}$ \textit{the drop of water follows the side of the biggest umbrella met. Note that in }$\rm{b)}$ \textit{the longest umbrellas are also the ones that are the narrowest.}}\label{fig:umbrellas}
% \label{umbrellas}
%\end{center}
%\end{figure}

\begin{figure}
\begin{center}
\includegraphics[scale=0.3]{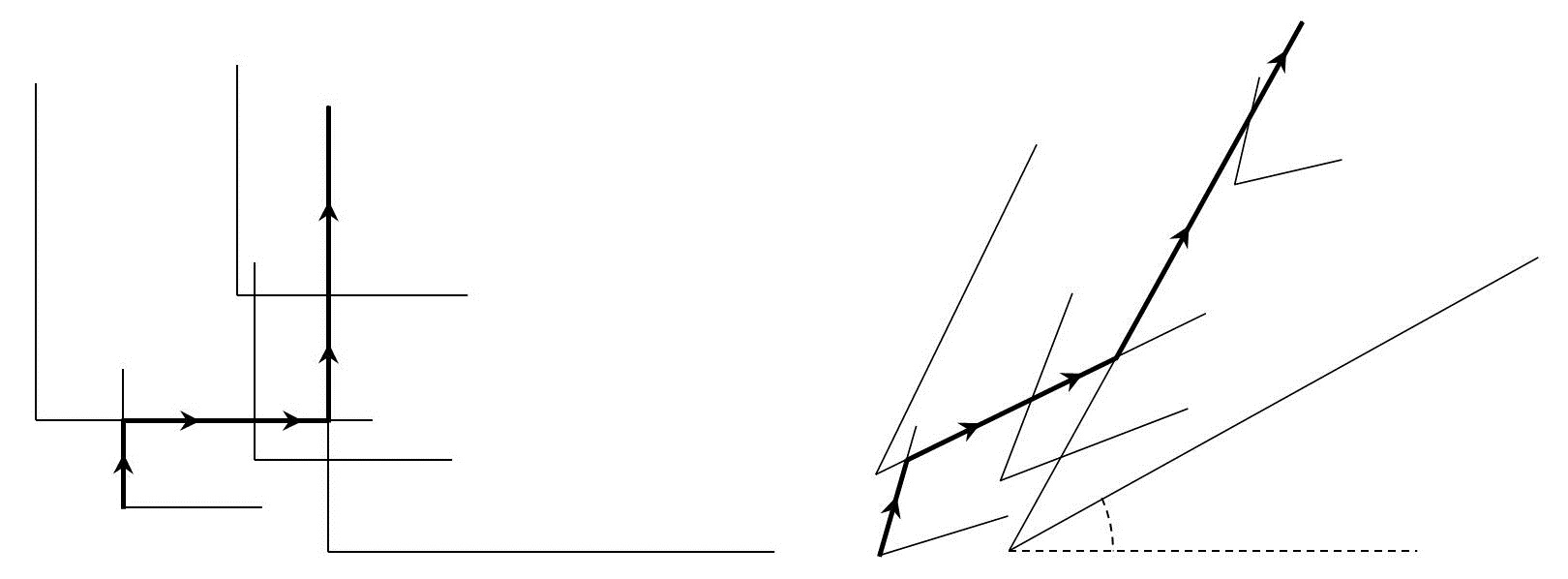}
 \caption{\textit{Both in the straight umbrellas case of Section \ref{sec:nospeed} and in the narrow umbrellas case of Section \ref{sec:speedpos}, the drop of water follows the side of the biggest umbrella met. Note that in the second case (right picture) the longest umbrellas are also the ones that are the narrowest.}}\label{fig:umbrellas}
 \label{umbrellas}
\end{center}
\end{figure}

We will now slightly modify the example seen in the previous subsection. Our aim is to build a new tree for which the behaviour of $h(0)$ is essentially the same as in the BZZ tree, but with a different shape of the graph. Roughly speaking, it will not allow to have long strips that are ''too horizontal'' or ''too vertical''. This feature and its importance will become more clear when we will describe the dynamics on these trees.

Fix suitable constants $\theta$ and $n_0\in\N$ such that $10 \leq\theta\leq n_0^2$ and so that following equation (\ref{eq:distribution}) makes sense. For every $x\in\Z^2$ consider i.i.d.~random variables $L(x)>1$ with atomless distributions fulfilling
\begin{equation}\label{eq:distribution}
  P(L(0)>t)=\frac{\theta\log t}{t^2}\qquad \mbox{for all } t\geq n_0.
\end{equation}

The new umbrellas we want to open 
%to repair the points of the lattice from the rain 
are a bit different from the tilde-umbrellas of the previous section. 

Define an umbrella of intensity $t$ as
\begin{equation}
 U_t=\bigcup_{i=1,2} U_{i,t}
\end{equation}
where
$U_{2,t}$ is the best $\Z^2$-grid lower approximation of the open segment of length $t$ that makes an angle of $\frac{\pi}{4}-\frac{1}{\log t}$ with the $x$-axis, living in the first quadrant and starting in the origin. $U_{1,t}$ is the reflection of $U_{2,t}$ with respect to the bisecting line of the first quadrant. % (see picture ???). 
$U_{1,t}$ and $U_{2,t}$ are the sides of the umbrella. Note that this time the intensity gives us the strength, the length but also the width of the umbrella. 
In particular, the longer the umbrella, the more narrow it is.

We can think once more that drops of rain pouring from every point of the lattice try to fall towards the up-right direction and that every time they reach a new vertex, they are deflected by the strongest umbrella that passes through that vertex (see Figure \ref{fig:umbrellas}).

In analogy with the straight-umbrellas case we define the strongest umbrella through $x$ perpendicular to direction $e_i$, for $i,j\in\{1,2\}$ and $i\not =j$, as
\begin{equation}\label{eqn:lambda}
 \lambda_i(x)=\sup_{y\in\Z^2:\,[x,x+e_j]\in\, y+U_{L(y)}}L(y).
\end{equation}

Note that since $L(0)>1$ and since we are taking the lower (for the first component) and upper (for the second) approximations of the segments described above, $[x,x+e_1]\in U_{2,L(x-e_1)}$ and $[x,x+e_2]\in U_{1,L(x-e_2)}$, so that the $\sup$ on the right hand side of \eqref{eqn:lambda} is taken over a non-empty set.
It requires slightly more work compared to the straight-umbrellas case to prove that it is also a.s. finite and therefore well defined.

We need some more notations. Similarly to \cite{BZZ}, for $m,n\in\Z$ call $S_m^n$ the slab 
$$
S_m^n=\big\{x=(x_1,x_2)\in\Z^2: m\leq x_1+x_2 \leq n \big\}.
$$

\begin{figure}
  \centering
  \setlength{\unitlength}{0.12\textwidth}  
  \begin{picture}(5,5)(0,1)
    \put(-2,0){\includegraphics[scale=0.50]{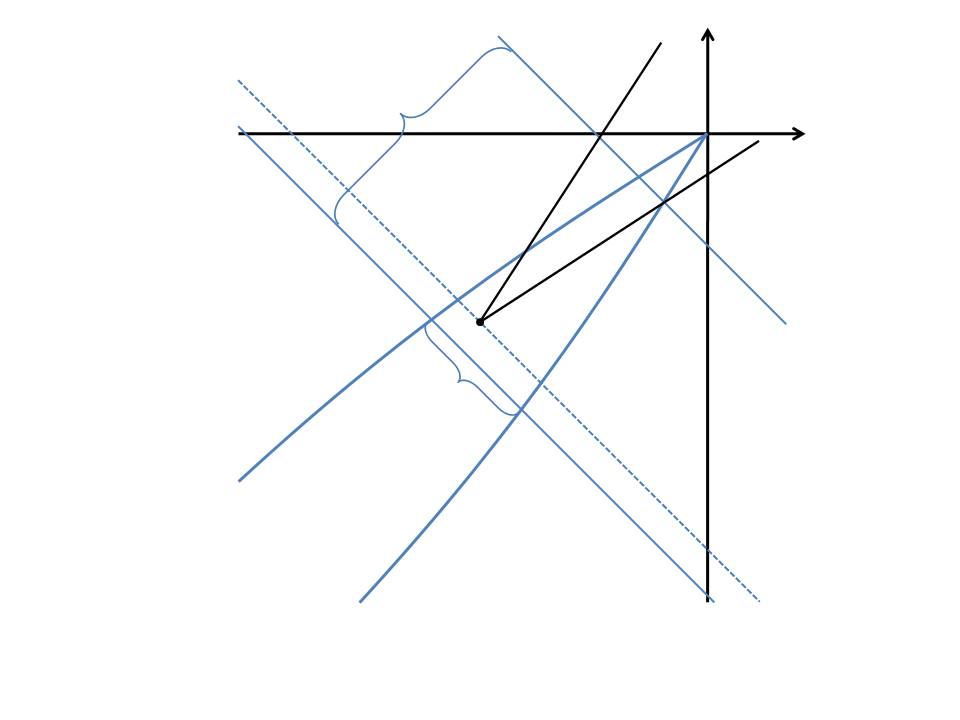}}
	 \put(1.4,2.6){$\frac{2n}{3\log n}$}
	 \put(0.89,5.095){$S_{-n}^{-n_0}$}
	 \put(0.1,5.3){$S^{-m}_{-m}$}
	 \put(1.83,3.2){$x$}
 	 \put(3.95,4.95){$O$}
	 \put(0.2,1.7){\Large{$G$}}
  \end{picture}
  \caption{\textit{The protecting area $G$ is the region of the plane from which we can have umbrellas that protect the origin. In particular, having a suitably strong umbrella starting in the part of $G$ delimited by the slab $S_{-n}^{-n_0}$ will ensure $h(0)<n$ with high probability.}}\label{gz}
\end{figure}
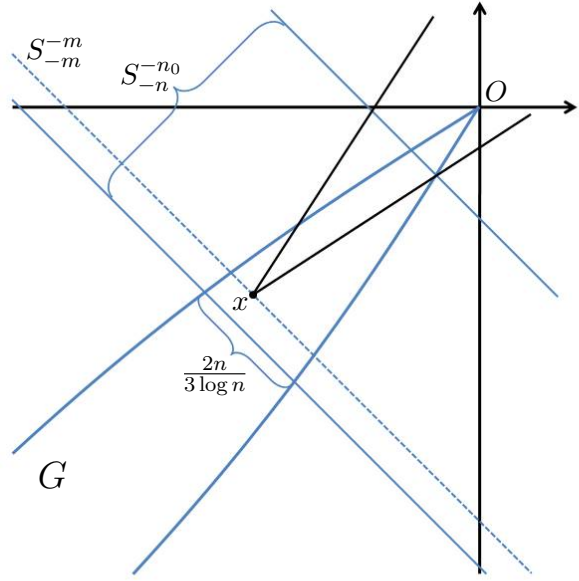

The protecting area $G$ (see Figure \ref{gz}) is defined as
\begin{align}
 G:=\Big\{&x=(x_1,x_2)\in-\N^2\big|\, \exists n\in\N: \nonumber\\
&x\in S_{-n}^{-n} \mbox{ and } -x_1\in\Big[y_n\cdot\cos\big(\tfrac \pi 4 -\alpha_n\big),y_n\cdot\cos\big(\tfrac \pi 4+\alpha_n\big)\Big]\Big\},
\end{align}
where $\alpha_n=\arctan\frac{\sqrt 2}{3\log n}$ and $y_n=\frac{n}{3\sqrt 2 \log n}\sqrt{2+9\log^2n}$. These %clumsy
values guarantee that every segment $S_{-m}^{-m}\cap G$ is $\frac{2m}{3\log m}$ long, and therefore contains $\frac{\sqrt 2 m}{3\log m}$ points of $\Z^2$ (up to one unit, at most).

%\begin{figure}
%\begin{center}
%\includegraphics[scale=0.50]{Good_zone}
% \caption{\textit{The protecting area $G$ is the region of the plane from which we can have umbrellas that protect the origin. In particular, having a suitably strong umbrella starting in the part of $G$ delimited by the slab $S_{-n}^{-n_0}$ will ensure $h(0)<n$ with high probability.}}\label{gz}
%\end{center}
%\end{figure}

Note that every umbrella $x+U_s$ with $x\in G$, $-(x_1+x_2)=n$ and $s\in[n,n^2]$, ``protects'' the origin $0$, meaning that $0$ lies inside the ``$\Z^2$-triangle'' generated by the sides $x+U_{1,s}$ and $x+U_{2,s}$.

\begin{lemma}
 There is a constant $c$ such that for $i=1,2$ and $t>n_0$,
\begin{equation}
 P(\lambda_i(0)>t)\leq c\frac{\log t}{t}.
\end{equation}
\end{lemma}
\begin{proof}
This is a straightforward calculation.
\begin{eqnarray*}
P(\lambda_i(0)>t) &\leq& 
C\int_t^\infty [s] \frac{\log s}{s^3} ds \\
\leq C \sum_{k=0}^\infty \int_{2^kt}^{2^{k+1}t} s \frac{\log s}{s^3} ds
&\leq& C \sum_{k=0}^\infty \frac{\log 2^kt}{2^k t}\\
= C \frac 1t\sum_{k=0}^\infty \frac{1}{2^k}[\log t + \log 2^k]
&\leq& c\,\frac{\log t}{t}.
\end{eqnarray*}
\end{proof}

Also in this case, the fact that the distributions of the $L(x)$'s are atomless guarantees the uniqueness of a direction $I(x)\in\{1,2\}$ such that
$$
\lambda_{I(x)}(x)=\min\{\lambda_i(x), i=1,2\}.
$$ 
For example, if $I(x)=1$, it means that the strongest vertical umbrella through $x$ is weaker than the strongest horizontal one. $I(x)$ is the direction which the drop of water will follow.

We can therefore define the new ancestral function
\begin{equation}\label{eqn:BSancestral}
 a(x)=x+e_{I(x)}.
\end{equation}
By its construction, it follows automatically that $a:\Z^2\to\Z^2$ is stationary and directed.

\begin{theorem}\label{thm:BS_tree}
 The random ancestral function described in \eqref{eqn:BSancestral} is such that
\begin{equation}\label{eqn:upper_bound_h_log}
 \limsup_{n\to\infty}\frac{n}{{\log^2 n}}P(h(0)>n)<\infty.
\end{equation}
\end{theorem}

Note that, even if \eqref{eqn:upper_bound_h_log} gives a slightly worse decay than \eqref{eqn:upper_bound_h}, the logarithmic correction will not affect the behaviour of the $\alpha$-logmoments of the conductances built on the different trees (cf.~Proposition \ref{prop:alpha_moments}).

\subsection{Proof of Theorem \ref{thm:BS_tree}}

 We closely follow the proof of Theorem 2 in \cite{BZZ}. 

We say that an umbrella $U$ penetrates a weaker umbrella $V$ in point $x\in\Z^2$ if one side of $U$ intersects one side of $V$ and $x$ is the upper-right point of their intersection. The following lemma bounds the probability that an umbrella of intensity $t$ starting in the origin gets penetrated by another stronger umbrella in a given point $z$.
\begin{lemma}\label{lemma:weak_umbrella}
Fix any $t>n_0$. Let $z\in\Z^2$ such that $[z,z+e_i]\in U_{j,t}$, for some $i,j\in\{1,2\}$. Then there exists a constant $c>0$ independent of $t$ such that
\begin{equation}
 P\big(I(z)\not=i|L(0)=t\big) \leq c\,\frac{\log{ t }}{ t }.
\end{equation}
\end{lemma}

\begin{proof}
For convenience, we shift the umbrella so that $z$ is translated to the origin. We look first at the event $E_k$ that the umbrella gets penetrated in the origin by an umbrella of intensity $s\in[k,k+1]$, for $k+1>t$. Note that such a penetrating umbrella can come only from $S_{-k-1}^{-1}$. Furthermore, on every $S_m^{m}$, $m\in \{-k-1,...,-1\}$, there are almost surely at most four  points that can generate it, since for all the others the slope of the sides would prevent them from penetrating the original umbrella in the origin.
Hence
\begin{align*}
P(E_k)&\leq \sum_{m=1}^{k+1}4\biggl(\frac{\theta\log k}{k^2}-\frac{\theta\log (k+1)}{(k+1)^2}\biggr)\leq c'\,\frac{\log k}{k^2},
\end{align*}
for some constant $c'$. 
It is now easy to see that
\begin{align*}
P\big(I(0)\not=i|L(-z)=t\big)\leq \sum_{k=\lfloor t\rfloor}^{\infty}\P(E_k)\leq c\,\frac{\log{t }}{t }.
\end{align*}  
\end{proof}
We now come to the proof of Theorem \ref{thm:BS_tree}.
For $n\geq n_0$, define now the random variables $M_n\in\{n_0-1,...,n\}$ as following:
\begin{equation}
M_n:=\max\big\{m\in\{n_0,...,n\}:\,\exists x\in S_{-m}^{-m}\cap G \mbox{ with }  m< L(x) <m^2\big\},
\end{equation}
with the convention $M_n=n_0-1$ whenever the set on the right hand side is empty.

Proving that, for some constant $c$,
\begin{equation}\label{eqn:em_en_estimates}
P(h(0)>m,\,M_n=m)\leq c\frac{\log^2 n}{n^2} \qquad \forall m=n_0,...,n,
\end{equation}
would imply
\begin{align}
P(h(0)>n)\leq\sum_{m=n_0-1}^n P(h(0)>m,\,M_n=m)\leq c \frac{\log^2 n}{n},
\end{align}
that is the statement of the theorem.

We first prove (\ref{eqn:em_en_estimates}) in the easy case $m=n_0-1$. 
\begin{align}
P(h(0)>n_0-1,&\,M_n=n_0-1)\leq P(M_n=n_0-1) \nonumber\\
  &=P\big(\,\mbox{for all } m=n_0,...,n \mbox{ and } x\in S_{-m}^{-m}\cap G,\,L(x)\not\in (m,m^2)\big) \nonumber\\
    &=\prod_{m=n_0}^n \big(1-P(L(0)\geq m)+P(L(0) > m^2)\big)^{\# (S_{-m}^{-m}\cap G) } \nonumber\\
      &=\prod_{m=n_0}^n {\Big(1-\frac{\theta\log m}{m^2}+\frac{\theta\log (m^2)}{m^4}\Big)}^{\frac{\sqrt2}{3}\frac{m}{\log m}} \nonumber\\
        &\leq\prod_{m=n_0}^n \Big(1-\theta\big(1-\tfrac 2{n_0^2}\big) \frac{\log m}{m^2}\Big)^{\frac{\sqrt2}{3}\frac{m}{\log m}}  \nonumber\\
          &\leq \e^{-\theta\big(1-\tfrac 2{n_0^2}\big)\frac{\sqrt 2}{3} \sum_{m=n_0}^n \frac{1}{m}} \leq c\, n^{-2}
\end{align}
by the choice of $\theta$, for some $c>0$.

For the more complex cases $m=n_0,...,n$ we faithfully follow \cite{BZZ} once again. For $i,j,r\in \Z$, $i\leq j$ and $x\in\Z^2$, define the events
\begin{equation}
A_i^j(x,r)=\big\{L(y)\not\in \big(-y\cdot \vec 1+r,(-y\cdot \vec 1+r)^2\big) \mbox{ for all } y\in S_i^j\cap (x+G)\big\}.
\end{equation}

%Our first task is to prove
%\begin{align}\label{eqn:first_part}
% \P(h(0)&>m,\,M_n=m) \nonumber\\
%&\leq \sum_{i=1,2}\E\Big[\#\{x\in S_m^{m}\cap -G:\,h(x)>m\};\,L(0)\in(m,m^2);\,A_{m-n}^{-1}\Big(\Big\lfloor\frac m 2\Big\rfloor e_i,m\Big)\Big],
%\end{align}

Firstly note that
\begin{align}
 P\big(h(0)>m,&\,M_n=m\big)%&= \nonumber\\
%&P(h(0)>m,\,L(y)\not\in(-y\cdot \vec 1+r,(-y\cdot \vec 1+r)^2) \; \forall y\in S_{-n}^{-m-1}\cap G,\,\exists x\in S_{-m}^{-m}\cap G:\,L(x)\in(m,m^2)) \nonumber\\
\leq \sum_{x\in S_{-m}^{-m}\cap G}P\big(h(0)>m,\,L(x)\in(m,m^2),A_{-n}^{-m-1}(0,0)\big)\nonumber\\
&=\sum_{x\in S_{-m}^{-m}\cap G}P\big(h(-x)>m,\,L(0)\in(m,m^2),A_{m-n}^{-1}(-x,m)\big)\nonumber\\
&=\sum_{x\in S_{m}^{m}\cap -G}P\big(h(x)>m,\,L(0)\in(m,m^2),A_{m-n}^{-1}(x,m)\big),\nonumber
\end{align}
where we have used stationarity to obtain the second line and we write $-G=\big\{x=(x_1,x_2):\,(-x_1,-x_2)\in G\big\}$.

\begin{figure}
  \centering
  \setlength{\unitlength}{0.12\textwidth}  
  \begin{picture}(6,5.8)
    \put(-0.2,-0.3){\includegraphics[scale=0.55]{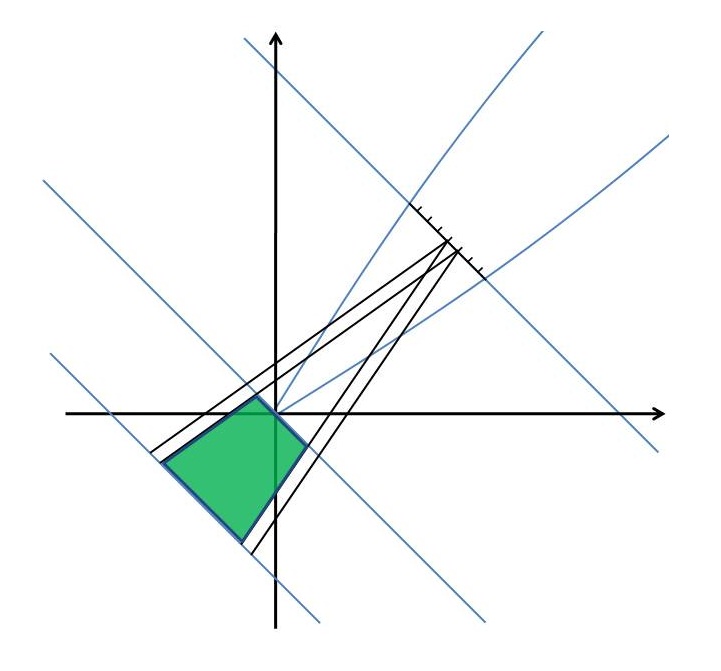}}
	 \put(4.3,0.2){$S_{0}^{0}$}
	 \put(5.55,1.45){$S^{m}_{m}$}
	 \put(2.75,0.2){$S^{-n+m}_{-n+m}$}
	 \put(3.9,3.7){$\hat x^j$}
	 \put(4.05,3.55){$\check x^j$}
	 \put(2.77,3.68){$\hat C^j$}
	 \put(2.97,3.63){\vector(1,-1){0.4}}
	 \put(4,2.4){$\check C^j$} 
	 \put(4,2.6){\vector(-1,1){0.4}}
	 \put(3.85,4.17){{$I_1$}}
	 \put(3.84,4.23){\vector(-1,-1){0.23}}
	 \put(5,5){{$-G$}}
  \end{picture}
  \caption{\textit{The event $E_j(m,n)$, involving only the points in $\hat C^j\cap\check C^j\cap S_{-n+m}^{-1}$ (the green area in the picture), is contained in $A_{m-n}^{-1}(x,m)$ for all $x\in I_j$.}}\label{pic:eight}
\end{figure}
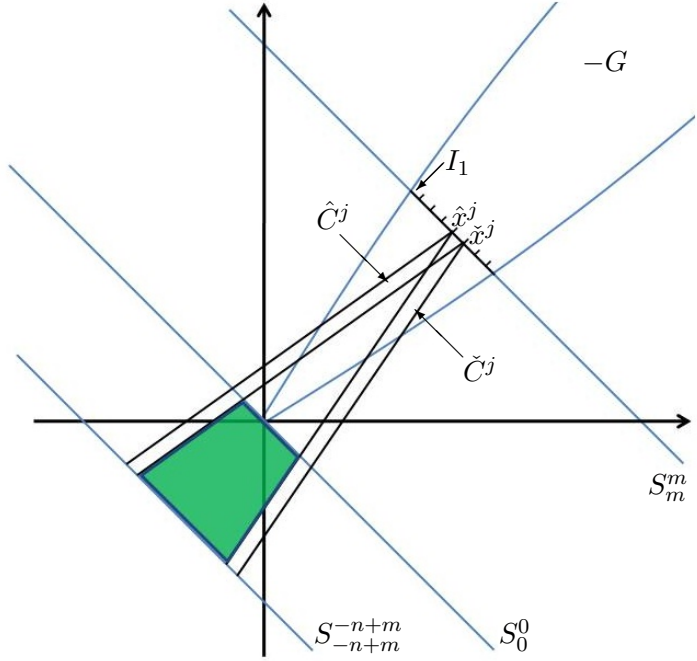

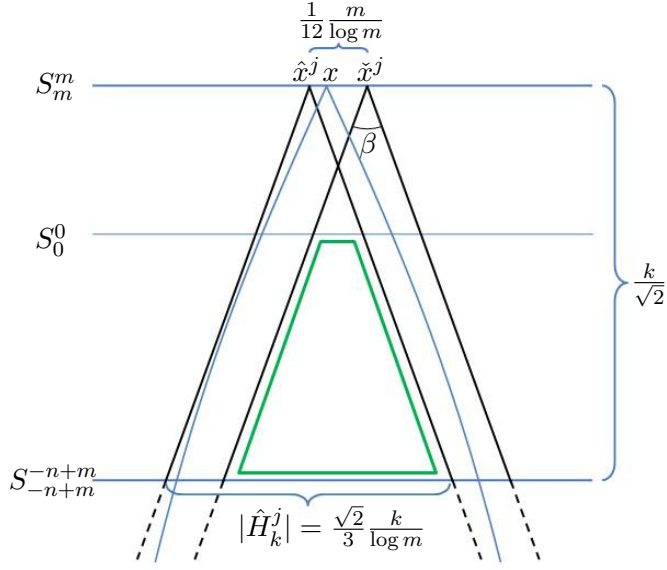
\begin{figure}
  \centering
  \setlength{\unitlength}{0.12\textwidth}  
  \begin{picture}(4,5)
    \put(-1.7,0){\includegraphics[scale=0.55]{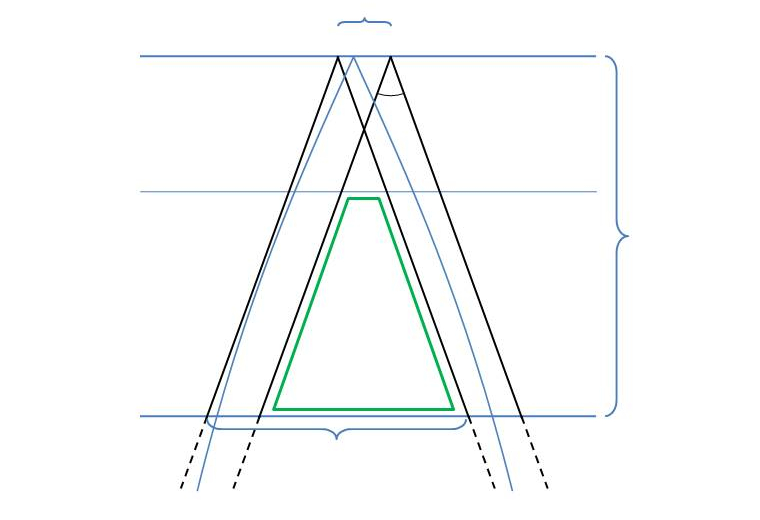}}
	 \put(-0.9,2.9){$S_{0}^{0}$}
	 \put(-0.9,4.2){$S^{m}_{m}$}
	 \put(-1.12,0.88){$S^{-n+m}_{-n+m}$}
	 \put(1.26,4.3){$\hat x^j$}
	 \put(1.8,4.3){$\check x^j$}
	 \put(1.51,4.3){$x$}
	 \put(1.3,4.75){$\tfrac{1}{12}\tfrac{m}{\log m}$}
	 \put(4.1,2.55){$\tfrac{k}{\sqrt 2}$} 
	 \put(0.8,0.48){{$|\hat H_k^j|=\tfrac{\sqrt 2}{3}\tfrac{k}{\log m}$}}
	 \put(1.82,3.7){$\beta$}
  \end{picture}
  \caption{\textit{The green area $\hat C^j\cap\check C^j\cap S_{-n+m}^{-1}$ is contained in $S_{-n+m}^{-1}\cap (x+G)$ for all $x\in I_j$ and it contains 'enough' points.}}\label{pic:pyramids}
\end{figure}

Consider now the segment joining the points in $S_{m}^{m}\cap -G$, divide it in eight parts of the same length (approximately $\frac 1 8 \frac{2m}{3\log m}$ long) and call them $I_1,...,I_8$ (see FIGURE \ref{pic:eight} and \ref{pic:pyramids}). For every $j\in\{1,...,8\}$, consider $\hat x^j$ and $\check x^j$, the points with respectively the highest and the lowest $y$-coordinate on $I_j$. Draw the infinite cones $\hat C^j$ and $\check C^j$ with angle $\beta=\arctan\big(\frac 2 {3\log m}\big)$ whose bisector makes an angle of $\frac {5}{4} \pi $ with the $x$-axis and with vertices $\hat x^j$ and $\check x^j$ respectively. Observe that the points in the area $\hat C^j\cap\check C^j\cap S_{-n+m}^{-1}$ are contained in $S_{-n+m}^{-1}\cap (x+G)$ for every $x\in I^j$. 
Therefore the event 
\begin{align*}
E_j(m,n):=\big\{L(y)\not\in \big(-y \cdot \vec 1 +m &,(-y\cdot \vec 1+m)^2\big) \\
&\mbox{ for all } y \in \hat C^j\cap\check C^j\cap S_{-n+m}^{-1}\big\} 
\end{align*}
is contained in the event $A_{m-n}^{-1}(x,m)$ for all $x\in I_j$.
Hence
\begin{align}\label{eqn:first_part}
 \sum_{x\in  S_{m}^{m}\cap -G}&P\big(h(x)>m,\,L(0)\in(m,m^2),A_{m-n}^{-1}(x,m)\big)\nonumber\\
&\leq
\sum_{j=1}^8 \sum_{x\in I_j}P(h(x)>m,\,L(0)\in(m,m^2),A_{m-n}^{-1}(x,m))\nonumber\\
&\leq \sum_{j=1}^8 \sum_{x\in I_j}P(h(x)>m,\,L(0)\in(m,m^2),E_j(m,n))\nonumber\\
&= \sum_{j=1}^8 E\Big[\#\{x\in I_j:h(x)>m\};\,L(0)\in(m,m^2);\,E_j(m,n)\Big].
\end{align}

The interval $(m,m^2)$ can be divided in a finite number of disjoint subintervals such that the $\Z^2$ approximation of every umbrella with intensity in a given subinterval looks the same at least up to the first $m$ edges. More precisely, there exists $M\in\N$ and there exist $\{m_1=m< m_2<...<m_M=m^2\}$ such that, for any $k\in\{1,2,...,M\}$, $\forall h,l\in (m_k,m_{k+1})$, one has $U_{h}|_m=U_{l}|_m$, where $U_{h}|_m$ is the umbrella of intensity $h$ whose sides are restricted to the first $m$ edges (going from bottom-left towards up-right).
Therefore, we can rewrite (\ref{eqn:first_part}) as
\begin{equation}\label{eqn:first_part2}
\sum_{j=1}^8\sum_{l=1}^{M-1} E\Big[\#\{x\in I_j:h(x)>m\};\,L(0)\in(m_l,m_{l+1});\,E_j(m,n)\Big].
\end{equation}

For any point $x\in S_{m}^{m}\cap -G$ to have $h(x)>m$, there must be a branch coming out of $x$ that perforates the protecting umbrella generated by the origin (since $L(0)\in(m,m^2)$). That is, at least one point $z$ on $U_{L(0)}|_m$ must be penetrated by another umbrella. On the other hand, every penetrated $z$ can give rise to at most one of such $x$'s. Hence, for any $l=1,...,M$, given $L(0)\in(m_l,m_{l+1})$,
\begin{align}
 \#\{x\in S_m^m\cap-G:\,h(x)>m\}\leq \sum_{i=1,2}\sum_{[z,z+e_i]\in U_{L(0)}|_m}\1_{\{I(z)\not=i\}}.
\end{align}

Plugging this in (\ref{eqn:first_part2}) gives
\begin{align}
 P&\big(h(0)>m,\,M_n=m\big)\nonumber\\
&\leq \sum_{j=1}^8\sum_{l=1}^M\sum_{i=1,2}\sum_{[z,z+e_i]\in U_{L(0)}|_m}P\big(I(z)\not=i,\,L(0)\in(m_l,m_{l+1}),\,E_j(m,n)\big).
\nonumber
\end{align}
The intersection of the first two events inside the last probability is not independent of $E_j(m,n)$, but there is a negative correlation between them. We obtain therefore the upper bound
\begin{align*}
 P&\big(h(0)>m,\,M_n=m\big)\nonumber\\
&\leq \sum_{j=1}^8\sum_{l=1}^M\sum_{i=1,2}\sum_{[z,z+e_i]\in U_{L(0)}|_m}P\big(I(z)\not=i,\,L(0)\in(m_l,m_{l+1})\big)P\big(E_j(m,n)\big).
\end{align*}
We can now directly compute the right hand side of last expression. For $[z,z+e_i]\in U_{L(0)}|_m$ we have, by Lemma \ref{lemma:weak_umbrella},
\begin{align}
P\big(I(z)=i;\,L(0)\in(m_l,m_{l+1})\big)&=\int_{m_l}^{m_{l+1}}P\big(I(z)=i|L(0)=t\big)\Big(\frac{\rm d}{{\rm d}t} P(L(0)\leq t)\Big){\rm d}t \nonumber\\
&\leq \int_{m_l}^{m_{l+1}} c\,\frac{\log{ t }}{ t } \frac{\theta}{t^3}(2\log t-1){\rm d}t\nonumber\\&
\leq K \frac{\log^2 m}{m^4}(m_{l+1}-m_l),
\end{align}
%\begin{align}
%\P\big(I(z)=i;\,L(0)\in(m,m^2)\big)&=\int_m^{m^2}\P\big(I(z)=i|L(0)=t\big)\Big(\frac{\rm d}{{\rm d}t}{\rm P}(L(0)\leq t)\Big){\rm d}t \nonumber\\
%&\leq \int_m^\infty c\,\frac{\log{ t }}{ t } \frac{\theta}{t^3}(2\log t-1){\rm d}t
%\leq K \frac{\log^2 m}{m^3},
%\end{align}theorem
for some constant $K$.

Summing over the directions $i=1,2$ and over all the $z\in\Z^2$ such that $[z,z+e_i]\in U_{L(0)}|_m$ and then summing over $l=1,...,M$, one is left with a factor of order $\frac{\log^2 m}{m^2}$.

In order to evaluate the probability of any $E_j(m,n)$, note that, for $k\geq m$, every $\hat C^j\cap\check C^j\cap S_{-k+m}^{-k+m}$ contains more than $\frac 1 5 \frac{k}{\log k}$ points of the lattice. In fact (see FIGURE \ref{pic:pyramids}),
%(see picture ????),
each cone $\hat C^j$ and $\check C^j$ intersects $S(k)$, the hyperplane containing $S_{-k+m}^{-k+m}$, on the segments $\hat H^j_k$ and $\check H^j_k$, each of length bigger than $\frac{k}{\sqrt 2}\cdot\frac{2}{3\log m}$ (they are, in fact, the double of the cathetus of a right triangle, whose opposite angle measures $\frac{\beta}{2}$ radians and with the other cathetus $\frac{k}{\sqrt 2}$ long). Since $\hat x^j$ and $\check x^j$ are roughly $\frac 18\frac{2m}{3\log m}$ far apart, the intersection of $\hat H^j_k$ and $\check H^j_k$ is longer than $\frac{\sqrt 2\,k}{3 \log m}-\frac 18\frac{2m}{3\log m}\geq \frac 1 3\frac{k}{\log k}$. Being the distance between close points on $\hat C^j\cap\check C^j\cap S_{-k}^{-k}$ equivalent to $\sqrt 2$, the total number of points is bigger than $\frac{1}{\sqrt 2}\frac 1 3\frac{k}{\log k}\geq \frac 1 5 \frac{k}{\log k}$.
 By the independence of the $(L(x))_{x\in\Z^2}$
\begin{align}
 P\big(E_j(m,n)\big)&= \prod_{k=1}^{n-m} \Big(1-P\big(L(0)\not\in (m+k, (m+k)^2)\big)\Big)^{\# (S_{-k}^{-k}\cap \hat C^j\cap\check C^j)}\nonumber\\
&\leq \prod_{k=m+1}^n \Big(1-\frac{\theta\log k}{k^2}\Big) ^{\frac 1 5\frac{k}{\log k}}\nonumber\\
&\leq \exp\Big\{-\frac{\theta}{5}{\sum_{k=m+1}^n}\frac 1k\Big\}\nonumber\\
&\leq \exp\Big\{-\frac{\theta}{5}{\int_{m+1}^n}\frac 1s{\rm d}s\Big\}\nonumber\\
&=\Big(\frac n {m+1}\Big)^{-\frac{\theta}{5}}.
\end{align}

Putting all together and reminding that $\theta\geq 10$, we finally obtain, for some constant $c$,
\begin{align}
 P\big(h(0)>m,\,M_n=m\big)&\leq c\frac{\log^2 m}{m^2}\Big(\frac n {m+1}\Big)^{-\frac{\theta}{5}}\nonumber\\
&\leq c (m+1)^{\frac \theta 5-2} n^{-\frac\theta 5}\log^2 m\nonumber\\
&\leq c\, n^{-2}\log^2n.
\end{align}

\section{The environment}\label{sec:env}

The two random trees constructed in the previous sections will provide, in some sense, the support for our random environments. In both cases, the $\omega$'s are constructed in the following way.

Sample a realization of the tree as described above. For every $z\in\Z^2$, the edge $\{z,a(z)\}$ will have a conductance value of $\omega_{\{z,a(z)\}}=\e^{(h(z)+1)^A}$, where $a:\Z^2\rightarrow \Z^2$ is the ancestral function used for constructing the sampled tree and $A>1$ is a constant. We set all the other conductances to be equal to one.

%We call a point $x\in\Z^2$ a \textit{leaf} of the tree if there is no $y\in\Z^2$ such that $a(y)=x$, where $a:\Z^2\rightarrow \Z^2$ is the ancestral function used for constructing the sampled tree. Of course, only one of the edges coming out of a leaf $x$ belongs to the tree, that is the edge $\{x,a(x)\}$, and we assign to it a conductance value of $\e$. Following the ancestral function we try to give to the next bond, for some $A>1$, a value $\omega(a(x),a(a(x)))=\e^{2^A}$, a value of $\e^{3^A}$ to the following one and so on. Whenever we try to assign two or more values to the same bond, we choose the bigger one. Therefore, for every $z\in\Z^2$, the edge $\{z,a(z)\}$ will have a value $\e^{(h(z)+1)^A}$.

For both the BZZ and the Diagonal tree, the conductances have infinite $\alpha$-logmoments for any $\alpha>1$. On the other hand, choosing appropriately the constant $A>1$, we can obtain conductances with finite $\alpha$-logmoments for $\alpha$ arbitrarily close to $1$ from below.

\begin{prop}\label{prop:alpha_moments}
 Take $\bar\alpha<1$. Then, the conductances of the random environments described above with $1<A< \tfrac 1 {\bar\alpha}$ are such that
\begin{equation}
 E\big[\log^\alpha \omega_e\big]<\infty \qquad\forall\alpha\leq\bar\alpha
\end{equation}
 and
\begin{equation}
 E\big[\log^\alpha \omega_e\big]=\infty \qquad\forall\alpha\geq 1.
\end{equation}
\end{prop}
\begin{proof}
We first prove it for the random environment built on the BZZ-tree support.
 \begin{align}\label{eqn:logmoments_series}
  E\big[\log^\alpha \omega_e\big]&=\int_0^\infty P(\log^\alpha \omega_e>t){\rm d}t\nonumber\\
				  &=\sum_{k=0}^\infty\int_{k^{\alpha A}}^{(k+1)^{\alpha A}} P(\log^\alpha \omega_e>t){\rm d}t\nonumber\\
				&\leq \sum_{k=0}^\infty P(\log^\alpha \omega_e>k^{\alpha A})((k+1)^{\alpha A}-k^{\alpha A})\nonumber\\
				&=\sum_{k=0}^\infty P(h(0)>k-1)((k+1)^{\alpha A}-k^{\alpha A}).
\end{align}
%Note that from the mean value theorem applied to the function $f(x)=x^{\alpha A}$, we have, for some $k\leq\xi\leq k+1$,
%$\alpha Ak^{\alpha A-1}\leq f'(\xi)=(k+1)^{\alpha A}-k^{\alpha A}\leq \alpha A(k+1)^{\alpha A-1}$.
%Furthermore, by 
By equations (\ref{eqn:lower_bound_h}) and (\ref{eqn:upper_bound_h}) we know that for all sufficiently large $k\in\N$, say $k\geq K$,
\begin{equation}
\frac c k\leq P(h(0)>k-1)\leq \frac{c'}{k}, \nonumber
\end{equation}
while the mean value theorem guarantees that
$$
\alpha A k^{\alpha A-1}\leq(k+1)^{\alpha A}-k^{\alpha A}\leq\alpha A {(k+1)}^{\alpha A-1}.
$$
Therefore, on the one hand, taking $\alpha\leq\bar\alpha$,
\begin{align}
  E\big[\log^\alpha \omega_e\big]&\leq C+\sum_{k=K}^\infty \frac {c'} k \alpha A(k+1)^{\alpha A-1}<\infty,\nonumber
\end{align}
where $C>0$ is the finite contribution of the first $K$ terms of the sum. On the other hand, when $\alpha\geq1$, we obtain with a minor modification of (\ref{eqn:logmoments_series})
\begin{align}
  E\big[\log^\alpha \omega_e\big]&> \sum_{k=K}^\infty \frac c k \alpha Ak^{\alpha A-1}=\infty. \nonumber
\end{align}

Note that the very same proof is valid for the random environment built over the Diagonal tree structure, since the $\log^2$-correction in Theorem \ref{thm:BS_tree} doesn't change the behaviour of the series (\ref{eqn:logmoments_series}).
\end{proof}

\begin{prop}\label{prop:follow_tree}
 For almost every environment $\omega$ sampled from the constructions of the previous section, the random walk among the conductances $\omega$ will eventually follow the tree. This means that almost surely there exists $\bar n< \infty$ such that for all $n\geq \bar n$, if $X_n=x$ then $X_{n+1}=a(x)$, where $a:\Z^2\rightarrow\Z^2$ is the ancestral function used to construct the tree underlying the environment.
\end{prop}

\begin{proof}
The probability that, starting in a point $x\in\Z^2$, the random walk will follow the tree forever is, by the independece of the jumps, bigger than
\begin{equation}\label{eqn:follow_tree}
\prod_{k=1}^\infty \frac{\e^{k^A}}{2\e^{(k-1)^A}+\e^{k^A}+1}.
\end{equation}
It is easy, in fact, to get convinced that this is a very pessimistic estimate. It represents the case in which we start from a leaf of the tree (that is, a vertex that is ancestor of no other vertices) and where every time $\omega_{\{X_n,a(X_n)\}}$ is of order $k$ (that is, equal to $\e^{k^A}$), then the two edges under and at the left of $X_n$ are of order $k-1$.

Call $T_1,T_2,...$ the times in which the random walk doesn't go in the direction of the ancestral function. After each of these times, a new attempt to follow the tree is performed.
Therefore if we show that the product (\ref{eqn:follow_tree}) is a constant strictly bigger than zero, than the sum of the probabilities of succeeding in following the tree in one of the attempts is infinite. By the Borel-Cantelli lemma, this means that almost surely there will be a finite time from which we will always follow the tree.

We are left to show that (\ref{eqn:follow_tree}) is bigger than zero, or, equivalently, that its $\log$ is bigger than $-\infty$:
\begin{align}
 \log \Big(\prod_{k=1}^\infty \frac{\e^{k^A}}{2\e^{(k-1)^A}+\e^{k^A}+1}\Big)&=
-\sum_{k=1}^\infty \log\Big(1+2\frac{\e^{(k-1)^A}}{\e^{k^A}}+\frac{1}{\e^{k^A}}\Big) \nonumber\\
&>-\sum_{k=1}^\infty \big(2\e^{(k-1)^A-k^A} +\e^{-k^A}\big)\nonumber\\
&>-\sum_{k=1}^\infty \big(2\e^{-A(k-1)^{A-1}}+\e^{-k^A}\big)>-\infty,
\end{align}
where we have used the mean value theorem for the bound $k^A-(k-1)^A\geq A(k-1)^{A-1}$.

\end{proof}

\begin{prop}\label{prop:bzzspeed}
The random walk among random conductances with environment built on the BZZ tree, as described above, has almost surely no limiting speed.
\end{prop}

\begin{prop}\label{prop:diagspeed}
The random walk among random conductances with environment built on the diagonal tree, as described above, has almost surely a limiting speed which is not zero.
\end{prop}

\begin{proof}[Proof of Proposition \ref{prop:bzzspeed}]From Proposition \ref{prop:follow_tree} we know that with probability $1$ there exists a finite time $\bar n$ from which the random walk will use only edges pointing the right or up direction with respect to its current position. Without loss of generality we can think this time to be time $0$.
In order to study the limiting speed of the process, we have to go back to the underlying structure of the tree on which we have built the environment. Note that every time we move one step in the direction of the ancestral function, we find several new umbrellas perpendicular to the step we have taken and a new parallel one. If the strongest perpendicular umbrella is stronger than any other umbrella on the direction of the previous step, the branch of the tree changes orientation; otherwise, it will continue in the same direction as before.

The distribution of the length $\bar L$ of the strongest new perpendicular umbrella met at each step in the direction of the ancestral function is easy to calculate:
\begin{align}
P(\bar L>t)&=P\big(\exists j\in\N \mbox{ such that } \tL((0,-j))>\max\{t,j\}\big)\nonumber\\
&=1-\prod_{j=1}^{\lfloor t\rfloor}P(\tL((0,-j))\leq t)
\prod_{j=\lfloor t\rfloor+1}^\infty P(\tL((0,-j))\leq j)\nonumber\\
&=1-\Big(1-\frac {\tilde\theta}{t^2}\Big)^{\lfloor t\rfloor}
\prod_{j=\lfloor t\rfloor+1}^\infty \Big(1-\frac{\tilde \theta}{j^2}\Big)
\end{align}
so that, by straightforward calculations, for $t$ sufficiently large,
\begin{equation}\label{eq:new_length}
\frac {c'}{t}\leq P(\bar L>t) \leq \frac {c''}{t},
\end{equation}
for some $c',c''>0$.

Following the tree in the direction of the ancestral function and considering only the strongest umbrellas through each point, call \textit{rush} a sequence of intersecting umbrellas, each bigger of the previous one, that determines a part of the final tree (note that a rush can well be formed by only one umbrella).

We will now proceed as follows: First of all we will prove that any rush is formed only by a finite number of umbrellas. We will then divide the time into accurately chosen intervals according to the rushes that the random walk will meet. Via Borel-Cantelli kind of arguments, we will prove that the random walk will follow very strong umbrellas for a sufficiently long time in infinitely many of these time-intervals. This will give every time a positive contribution to the velocity up to that point, showing that the velocity oscillates infinitely often and therefore cannot have a limit.  

We want to show now that any rush is formed only by a finite number of umbrellas almost surely. 
In fact, any side of an umbrella of length $k>2n_0^2$ meets neither a stronger perpendicular nor a stronger parallel umbrella with probability bigger than
%\P(\mbox{For $\lfloor k\rfloor$ times the stronger perpendicular and the new parallel umbrellas are $<k$})>
\begin{align}\label{eq:leave_a_rush}
\Big(1-\frac{c''}{k}\Big)^k\Big(1-\frac{\tilde\theta}{k^2}\Big)^k>\e^{-2c''-1},
\end{align}
that is a constant strictly bigger than $0$ and independent of $k$.
Therefore, the probability of $Q:=\{ \mbox{a rush is formed by more than $N$ umbrellas}$ $\mbox{longer than $2n_0^2$}\}$ is bounded by

%Starting on any rush, the probability of leaving it (that is, of travelling the whole length of one of the umbrellas without meeting a stronger one) after having visited $N\in \N$ different umbrellas is

\begin{align}
P(Q)&<\big(1-c\e^{-2c''}\big)^N,
\end{align}
for some $c>0$.

This means that, in almost every realization of the BZZ tree, each branch is determined only by rushes of finitely many umbrellas. From the point of view of the random walk that from time $0$ on will deterministically follow the ancestral function, this means that the walker will leave any underlying rush in finite time for almost every tree. Given a realization of the walk, call $\tau(1)\in\N$ the time in which the random walk leaves the first rush, $\tau(2)$ the time in which it leaves the second one and so on. $\tau(1)<\tau(2)<...$ is a sequence of (almost surely finite) integer times that goes to infinity.

Fix $T>1$ and define the times $T_1=T$, $\tau_1=\min_{i=1,2,...}\{\tau(i):\, \tau(i)>T_1\}$ and recursively
\begin{align}
T_k&=\tau_{k-1}+\tau_{k-1}T^k &\forall k>1,\nonumber\\
\tau_k&=\min_{i=1,2,...}\{\tau(i):\, \tau(i)>T_{k-1}\} &\forall k>1.
\end{align}

Our aim is now to show that, in the intervals of the form $(T_{k-1},T_k)$, the longest umbrella met is of length of the order $\tau_{k-1}T^k$.
We do not want the longest umbrella to be much longer than this, otherwise it could ``interfere'' with the next intervals (that is, in order to simplify the forthcoming calculations we want that from some point on the longest umbrella met in an interval has not been already met in a previous interval): Consider the event 
\begin{align*}
E_k=\{ &\mbox{In the interval $(T_{k-1},T_k)$ the longest umbrella met} \\ 
&\mbox{ is stronger than $\tau_k T^{k+1}$}\}.
\end{align*}
Its probability can be bounded from above by
\begin{align*}
 \P(E_k)&<1-\Big(1-\frac {c''}{\tau_kT^{k+1}}\Big)^{T_k} \\
        &<1-\e^{-\frac{2c''}{T^{k+1}}}\\
	&<\frac{c}{T^{k+1}},
\end{align*}
for some constant $c>0$, since $T_k\leq \tau_k$. By the Borel-Cantelli lemma, $\P(E_k \mbox{ i.o.})=0$.

On the other hand, we don't want the longest umbrella to be shorter than that. This is because we want it to be long a positive fraction of the entire time interval $(T_{k-1},T_k)$. In fact, the interval $(T_{k-1},T_k)$ is longer than $\tau_{k-1}T^k$. Furthermore, we want the random walk to follow this umbrella for a positive fraction (say an $\varepsilon>0$ fraction) of its length before leaving the time interval. This two events guarantee a relevant contribution to the speed up to time $T_k$.
Therefore take, for a fixed $\varepsilon>0$ small,
\begin{align*}
 F_k&=\big\{ \mbox{In the interval $(T_{k-1},T_k(1-\varepsilon))$ the strongest umbrella met is stronger } \\
 &\mbox{than } \tau_{k-1}T^{k} \mbox{ and is stronger than the strongest umbrella in } (T_k(1-\varepsilon), T_k) \big\}.
\end{align*}

By the independence of the new umbrellas discovered at each step, we have, for all $k\in\N$,
\begin{align}
 \P(F_k)&>(1-\varepsilon) \P(\mbox{one of the $T_k-T_{k-1}$ umbrellas is longer than $\tau_{k-1}T^k$})\nonumber\\        
&=(1-\varepsilon)\Big(1-\Big(1-\frac {c'}{\tau_{k-1}T^{k}}\Big)^{T_k-T_{k-1}} \Big(1-\frac {\tilde\theta}{\tau_{k-1}^2T^{2k}}\Big)^{T_k-T_{k-1}} \Big)\nonumber \\
&>(1-\varepsilon)\Big(1-\Big(1-\frac {c'}{\tau_{k-1}T^{k}}\Big)^{T_k-T_{k-1}}\Big)\nonumber \\
&>(1-\varepsilon)\Big(1-\e^{-\frac{c'} 2\frac{1}{\tau_{k-1}T^k} (\tau_{k-1}+\tau_{k-1}T^k)}\Big)\nonumber\\
		 &=C
\end{align}
where $C>0$ is a constant not depending on $k$. By the second Borel-Cantelli lemma, there are almost surely infinitely many intervals $(T_{k-1},T_{k})$ for which $F_k$ happens.

 Hence, almost surely there exists a $\bar k\in \N$ (depending eventually on the realization of the environment and of the random walk) such that $E_k$ does not happen for every $k>\bar k$ while $F_k$ holds infinitely many times.
Take now the strongest umbrella met up to time $T_{\bar k}$. Its length $L>0$ is almost surely finite, so that $\kappa:=\min\{k:T_k>T_{\bar k}+L\}$ is well defined.
Note that $\forall k>\kappa+1$, in the interval $(T_k,T_{k+1})$ there is no umbrella longer than $T^{k+1}\tau_k$ met in the past.

Take the infinite subsequence $\kappa<k_1<k_2<...$ such that $F_{k_i}$ holds true for every $i\in\N$ and such that the longest umbrella met in the $k_{i}$'th interval $(T_{k_i-1},T_{k_i})$ is followed by the random walk at least for a positive fraction $0<\eta<\epsilon$ of its length. Note that since there is no longer umbrella coming from a previous interval, once the random walk meets this umbrella it follows it until its end or at least until the end of the interval itself, and the probability of meeting the umbrella before the last $\eta$ fraction of its length is strictly positive. This implies that we have such a sequence $(k_i)_{i=1,2,...}$ almost surely.

Suppose now that a limiting speed $v=(v[1],v[2])$ existed. We want to show that in each of those intervals there is at least one time $t$ at which the ratio $\sum_{j=1}^t X_j/t$ is far from $v$, bringing to a contradiction.
Call $t_i\in\N$ the time at which the longest umbrella of the interval $(T_{k_i-1},T_{k_i})$ is met and $t^i=\{ \mbox{``Time of the last point of the umbrella''}\wedge T_{k_i}\}$. By definition, this umbrella is longer than $\tau_{k_i-1}T^{k_i}$, it is met before time $T_{k_i}(1-\varepsilon)$ and before the last $\eta$-fraction of its length.
Call
$$
\frac {1}{t_i}\sum_{j=1}^{t_i}X_j=(v_{i}[1],v_{i}[2])=:v_{i}
$$
and
$$
\frac {1}{t^i}\sum_{j=1}^{t^i}X_j=(v^{i}[1],v^{i}[2])=:v^{i}
$$
the partial speeds up to time $t_i$ and $t^i$ respectively. Without loss of generality suppose that we met the longest umbrella on its horizontal side. Note that
$$
v_i[1]=\frac{1}{t_i}\big(v^i[1]t^i-t^i+t_i\big)
$$
and that 
$$
\frac{t^i}{t_i}>1+\frac{\eta \tau_{k_i-1}T^{k_i}}{(1-\varepsilon)\tau_{k_i-1}(T^{k_i}+1)}>1+\frac{\eta}{2(1+\varepsilon)}=:\beta>0.
$$

Further suppose $v[1],v[2]\not \in\{0, 1\}$. Then if $v[1]>v^i[1]$
\begin{align}
|v[1]-v_i[1]|=v[1]-v^i[1]\tfrac{t^i}{t_i}+\tfrac{t^i}{t_i}-1>(\beta-1)(1-v[1])>0,
\end{align}
while if $v[1]\leq v^i[1]$
\begin{align}
 \max\big\{|v[1]-v^i[1]|,|v[1]-v_i[1]|\big\}&\geq \tfrac{1}{2}(v^i[1]-v_i[1])\nonumber\\
	&=\tfrac{1}{2}\big(v^i[1]-v^i[1]\tfrac{t^i}{t_i}+\tfrac{t^i}{t_i}-1\big)\nonumber\\
	&>\tfrac{1}{2}(\beta-1)(1-v[1]).
\end{align}
In both cases the distance from the limiting speed is bigger than a constant that is independent of $k_i$ and strictly bigger than zero.

The cases $v[1]=1$ and $v[1]=0$ have probability $0$.
In fact, the probability of meeting in any interval $(T_{k_{i-1}},(1-\varepsilon )T_{k_i})$ a vertical (respectively, horizontal) umbrella of order $\tau_{k_i-1}T^{k_i}$ that is stronger of any other horizontal (vertical) umbrellas met before (and of following it for a time of $O(t)$) is strictly positive, for the reasons mentioned above.

\end{proof}

\begin{proof}[Proof of Proposition \ref{prop:diagspeed}]
Let $v=(0.5,0.5)$. We claim that, almost surely, 
\[
\lim_{n\to\infty}\frac{X_n}{n}=v.
\]
As in the previous proof, let $\bar n$ be such that for every $n>\bar n$, we have $X_{n+1}=a(X_n)$, where $a$ is the ancestral function. By Proposition \ref{prop:follow_tree}, we know that $\bar n$ is almost surely finite also in the present case. 
We need to prove that for every $\varepsilon>0$ there exists a (random) finite $M$ such that for every $n>M$, we have $\|X_n/n-v\|<\varepsilon$, where we write $\|\cdot\|$ for, e.g., the usual $1$-norm. To this end, we need to understand the various umbrellas that the random walk traverses. 
So, let $\varepsilon>0$. By the construction of the diagonal tree, there exists $K>0$ such that for every umbrella which is stronger than $K$, for every two points $x$ and $y$ on the umbrella whose distance is larger than some $U=U(\varepsilon)$, we have 
\begin{equation}\label{eq:slp}
\left\|\frac{y-x}{\|y-x\|}-v\right\|<\varepsilon.
\end{equation}
%For umbrellas which are not stronger than $K$, their distribution is symmetric w.r.t. the diagonal, and their directions are i.i.d.~and therefore they give an average of $v$.

Let $\alpha(n)=\langle X_n,(-1,1) \rangle$ be the (signed) distance of $X_n$ from the diagonal. We will prove that almost surely,
\begin{equation}\label{eq:supbound}
\limsup_{n\to\infty} \frac{\alpha(n)}{n}\leq \varepsilon,
\end{equation}
and similarly
\begin{equation}\label{eq:infbound}
\liminf_{n\to\infty} \frac{\alpha(n)}{n}\geq -\varepsilon.
\end{equation}

To see \eqref{eq:supbound}, let $n_1:=\inf\{n:\alpha(n)>\varepsilon n\}\leq\infty$, and let 
\[
n_{k+1}:=\inf\{n>n_k+U(\epsilon): \alpha(n)-\varepsilon n> \alpha(n_k)-\varepsilon n_k\}\leq \infty.
\]
Next we prove that almost surely
\begin{equation}\label{eq:n_k_limit}
\lim_{k\to\infty} \frac{n_k}{k} = \infty.
\end{equation}

Indeed, fix  $L>K$. Let $E_k=\{n_{k+1}-n_k>\tfrac L 2\}$.  We claim that 
\begin{equation}\label{eq:last_equation}
P(E_k |\, X_1,...,X_{n_k})\geq \frac{C\log L}{L}
\end{equation}
for some constant $C>0$. In fact, $E_k$ happens if the point $X_{n_k}$ is on the first half of the side of an umbrella of strength larger than $L$, and the walker follows this umbrella to its end. Then, in order to justify \eqref{eq:last_equation}, it is enough  to observe that $X_{n_k}$ is exposed to the lower side of umbrellas that are independent from the past of the walk (meaning that we do not have information on their distribution) and to compute \eqref{eq:new_length} and \eqref{eq:leave_a_rush} exactly as before for the case of the diagonal tree. 
Note that the sequence $Y_k=\1_{E_k}$ dominates a sequence of i.i.d. Bernoulli random variables with mean $\tfrac{C\log L}{L}$, and hence, almost surely,
$$
\liminf_{k\to\infty}\frac{1}{k}\sum_{j=1}^{k}Y_k\geq \frac{C\log L}{L}.
$$
This, together with $n_k\geq\tfrac L 2\sum_{j=1}^k Y_k$ and the arbitrariness of $L$ gives \eqref{eq:n_k_limit}.

\eqref{eq:supbound} now follows easily: note that $\max_{n\leq n_k}\alpha(n)-\varepsilon n \leq k U(\varepsilon)$. Then, with $k_n = \max\{k:\, n_k\leq n\}$, we have
\begin{align*}
\limsup_{n\to\infty} \frac{\alpha(n)-\varepsilon n}{n}
&\leq \limsup_{n\to\infty} \frac{1}{n}\big(U(\varepsilon)(k_n+1)\big)\\
&=U(\varepsilon)\limsup_{n\to\infty}\frac{k_n}{n} =0
\end{align*}
since $\tfrac{n_k}{k}\to \infty$ implies $\tfrac{k_n}{n}\to 0$.
 To see \eqref{eq:infbound} note that the entire system is invariant to reflection $(x,y)\to(y,x)$.

\end{proof}

\section{Acknowledgment}

The authors would like to thank the referee for a nice and very fast review.
The second author would also like to thank Simon Wasserroth from TU Berlin for useful discussions on worms.
The research of the first author was partially supported by the Humboldt foundation and by ERC grant StG 239990.

\bibliography{speed}
\bibliographystyle{alpha}

\end{document}